\documentclass[12pt]{amsart}
\usepackage{graphicx}
\usepackage{amsmath,amsthm,amssymb,enumerate}
\usepackage{euscript,mathrsfs}
\usepackage[citecolor=blue,colorlinks=true]{hyperref}
\usepackage[left=3cm,right=3cm,top=4.5cm,bottom=4.5cm]{geometry}
\usepackage{color}
\catcode`\@=11 \@addtoreset{equation}{section}

\catcode`\@=12

\allowdisplaybreaks

\newtheorem{Theorem}{Theorem}[section]
\newtheorem{Proposition}[Theorem]{Proposition}
\newtheorem{Lemma}[Theorem]{Lemma}
\newtheorem{Corollary}[Theorem]{Corollary}

\theoremstyle{definition}
\newtheorem{Definition}[Theorem]{Definition}

\newtheorem{Remark}[Theorem]{Remark}

\newcommand{\bTheorem}[1]{
\begin{Theorem} \label{T#1} }
\newcommand{\eT}{\end{Theorem}}

\newcommand{\bProposition}[1]{
\begin{Proposition} \label{P#1}}
\newcommand{\eP}{\end{Proposition}}

\newcommand{\bLemma}[1]{
\begin{Lemma} \label{L#1} }
\newcommand{\eL}{\end{Lemma}}

\newcommand{\bCorollary}[1]{
\begin{Corollary} \label{C#1} }
\newcommand{\eC}{\end{Corollary}}

\newcommand{\bRemark}[1]{
\begin{Remark} \label{R#1} }
\newcommand{\eR}{\end{Remark}}

\newcommand{\bDefinition}[1]{
\begin{Definition} \label{D#1} }
\newcommand{\eD}{\end{Definition}}

\newcommand{\bfphi}{\boldsymbol{\varphi}}

\newcommand{\bFormula}[1]{
\begin{equation} \label{#1}}
\newcommand{\eF}{\end{equation}}

\newcommand{\Ov}[1]{\overline{#1}}

\newcommand{\aleq}{\lesssim}

\newcommand{\vr}{\varrho}

\newcommand{\vu}{\vc{u}}
\newcommand{\vm}{\vc{m}}

\newcommand{\vc}[1]{{\bf #1}}

\newcommand{\Div}{{\rm div}_x}
\newcommand{\Grad}{\nabla_x}

\newcommand{\dx}{\,{\rm d} {x}}

\newcommand{\dt}{\,{\rm d} t }

\newcommand{\vU}{\vc{U}}

\newcommand{\intO}[1]{\int_{\Omega} #1 \, \dx}

\newcommand{\D}{{\rm d}}

\definecolor{Cgrey}{rgb}{0.85,0.85,0.85}
\definecolor{Cblue}{rgb}{0.50,0.85,0.85}
\definecolor{Cred}{rgb}{1,0,0}
\definecolor{fancy}{rgb}{0.10,0.85,0.10}

\allowdisplaybreaks

\newcommand\Cbox[2]{%
    \newbox\contentbox%
    \newbox\bkgdbox%
    \setbox\contentbox\hbox to \hsize{%
        \vtop{
            \kern\columnsep
            \hbox to \hsize{%
                \kern\columnsep%
                \advance\hsize by -2\columnsep%
                \setlength{\textwidth}{\hsize}%
                \vbox{
                    \parskip=\baselineskip
                    \parindent=0bp
                    #2
                }%
                \kern\columnsep%
            }%
            \kern\columnsep%
        }%
    }%
    \setbox\bkgdbox\vbox{
        \color{#1}
        \hrule width  \wd\contentbox %
               height \ht\contentbox %
               depth  \dp\contentbox
        \color{black}
    }%
    \wd\bkgdbox=0bp%
    \vbox{\hbox to \hsize{\box\bkgdbox\box\contentbox}}%
    \vskip\baselineskip%
}

\date{}

\begin{document}


\title{Generalized solutions to models of inviscid fluids}

\author{Dominic Breit}
\address[D. Breit]{Department of Mathematics, Heriot-Watt University, Riccarton Edinburgh EH14 4AS, UK}
\email{d.breit@hw.ac.uk}

\author{Eduard Feireisl}
\address[E.Feireisl]{Institute of Mathematics AS CR, \v{Z}itn\'a 25, 115 67 Praha 1, Czech Republic
\and Institute of Mathematics, TU Berlin, Strasse des 17.Juni, Berlin, Germany }
\email{feireisl@math.cas.cz}
\thanks{The research of E.F. leading to these results has received funding from
the Czech Sciences Foundation (GA\v CR), Grant Agreement
18--05974S. The Institute of Mathematics of the Academy of Sciences of
the Czech Republic is supported by RVO:67985840.}

\author{Martina Hofmanov\'a}
\address[M. Hofmanov\'a]{Fakult\"at f\"ur Mathematik, Universit\"at Bielefeld, D-33501 Bielefeld, Germany}
\email{hofmanova@math.uni-bielefeld.de}
\thanks{M.H. gratefully acknowledges the financial support by the German Science Foundation DFG via the Collaborative Research Center SFB1283.}

\begin{abstract}

We discuss several approaches to generalized solutions of problems describing the motion of inviscid fluids. We propose a new concept of \emph{dissipative} solution to the compressible Euler system based on a careful analysis of possible oscillations and/or concentrations in the associated generating sequence. Unlike the conventional measure--valued solutions
or rather their expected values, the dissipative solutions comply with a natural compatibility condition -- they are classical solutions
as long as they enjoy certain degree of smoothness.

\end{abstract}


\date{\today}

\maketitle


\section{Introduction}
\label{I}

We consider a mathematical model of an \emph{inviscid compressible fluid} with the mass density $\vr = \vr(t,x)$ moving with the velocity $\vu = \vu(t,x)$. Thermal effects
being neglected, the evolution of the fluid is governed by the \emph{Euler system}:
\begin{equation} \label{i1}
\begin{split}
\partial_t \vr + \Div (\vr \vu) &= 0, \\
\partial_t (\vr \vu) + \Div (\vr \vu \otimes \vu) + \Grad p(\vr) &= 0.
\end{split}
\end{equation}
The quantity $p = p(\vr)$ is the pressure. We suppose the internal energy $e = e(\vr)$ is related to the pressure through the formula
\begin{equation} \label{i2}
(\gamma - 1) \vr e(\vr) = p(\vr),
\end{equation}
$\gamma > 1$ is the adiabatic constant. The total energy of the fluid is given by
\begin{equation} \label{i3}
E(\vr, \vu) = \varrho \left[ \frac{1}{2} |\vu|^2 + e(\vr) \right].
\end{equation}
If not otherwise stated, we suppose the fluid occupies a bounded domain $\Omega \subset R^d$, $d=2,3$ with impermeable boundary:
\begin{equation} \label{i4}
\vu \cdot \vc{n}|_{\partial \Omega} = 0.
\end{equation}
The initial state of the system is given:
\begin{equation} \label{i5}
\vr(0, \cdot) = \vr_0, \ \vr \vu(0, \cdot) = \vm_0.
\end{equation}

The Euler system has been derived from the basic principles of continuum mechanics on condition that all quantities appearing in
\eqref{i1} are at least continuously differentiable and the density $\vr$ is bounded below away from zero. If the initial data belong to this class then the Euler system admits local--in--time smooth solutions, see e.g. Tani \cite{TAN}. The life span of such solution, however, is finite for a fairly general class of the initial data, see Smoller \cite{SMO}.

To continue solutions globally in time, the concept of \emph{weak} solution is introduced, where all derivatives in \eqref{i1}
are understood in the sense of distributions. It is also more convenient to reformulate the problem in the \emph{conservative}
variables $\vr$ and $\vm = \vr \vu$:
\begin{equation} \label{i6}
\begin{split}
\partial_t \vr + \Div \vm &= 0, \\
\partial_t \vm + \Div \left( \frac{\vm \otimes \vm}{\vr} \right) + \Grad p(\vr) &= 0.
\end{split}
\end{equation}

Weak solutions are not unique unless a suitable admissibility criterion is imposed. In the context of the Euler system, it is customary to require the energy inequality
\begin{equation} \label{i7}
\partial_t E(\vr, \vm) + \Div \left[ \left( E(\vr, \vm) + p(\vr) \right) \frac{\vm}{\vr} \right] \leq 0,
\end{equation}
where
\begin{equation} \label{i8}
E(\vr, \vm) = \frac{1}{2} \frac{|\vm|^2}{\vr} + \vr e(\vr).
\end{equation}
In view of \eqref{i2}, we obtain
\[
p(\vr) = a \vr^\gamma, \ \vr e(\vr) \equiv P(\vr) = \frac{a}{\gamma - 1} \vr^\gamma.
\]
{Indeed the internal energy $\vr e$ coincides (modulo a linear function) with the pressure potential
$P = P(\vr)$}:
\[
\vr e(\vr) = P(\vr), \ \mbox{where the latter satisfies}\ P'(\vr) \vr - P(\vr) = p(\vr).
\]
In particular, the energy $E$ is a \emph{convex} function of $\vm$ and $\vr$.

Even if \eqref{i7} is imposed as an extra admissibility constraint, the Euler system remains ill--posed
at least for $N=2,3$. As a matter of fact, there exist
Lipschitz initial data for which \eqref{i6}, \eqref{i7} admits infinitely many solutions on a given time interval $(0,T)$,
see Chiodaroli et al. \cite{ChiDelKre}, \cite{ChiKre}, and \cite{ChKrMaSwI}. Still the question of the \emph{existence} of
global--in--time weak solutions
to \eqref{i6}, \eqref{i7} for \emph{general} initial data remains open.

Our goal is to present several concepts of generalized solutions to the Euler system and discuss their basic properties. In particular, we address the question of \emph{compactness} of the solution set and its \emph{stability} with respect to perturbations. Finally, we
introduce a new concept of \emph{dissipative} solution to the Euler system.

The paper is organized as follows. In Section \ref{O} we discuss the problem of compactness of the solution set of the compressible Euler system. Section \ref{M} presents a short review of various concepts of the so--called measure--valued solutions. In Section \ref{D}, we introduce a new concept of dissipative solutions. In Section \ref{C} we introduce a generating sequence and show existence of
a dissipative solution to the Euler system for a fairly general class of initial data. Various properties of dissipative solutions including weak--strong uniqueness and conditional regularity are discussed in Section \ref{P}. The paper is concluded in Section
\ref{S} by introducing
admissible dissipative solutions that maximize the mechanical energy dissipation.

\section{Oscillatory solutions}
\label{O}

As revealed by the method of convex integration, bounded sets of solutions to the Euler system may not be precompact even with respect to the natural weak topology, cf. e.g. De Lellis and Sz\' ekelyhidi \cite{DelSze3}. Indeed we claim the following result.

\begin{Proposition} \label{PO1}

Let $\Omega \subset R^N$, $N=2,3$ be a bounded domain. Let $\vr_0 \in L^\infty (\Omega)$, $\vr_0 > 0$ be given.

Then there exists a sequence of weak solutions $[\vr_n,\vm_n]$ to the Euler system \eqref{i6} in $(0,T) \times \Omega$ with $\vr_n=\vr_n(x)$ such that
\begin{equation} \label{O1}
\vr_n \to \vr_0 \ \mbox{weakly-(*) in}\ L^\infty(\Omega),\
\vm_n \to 0 \  \mbox{weakly-(*) in}\ L^\infty((0,T) \times \Omega; R^N),
\end{equation}
\begin{equation} \label{O2}
\liminf_{n \to \infty} \intO{ |\vr_n - \vr_0| } > 0.
\end{equation}

\end{Proposition}

\begin{Remark} \label{OR1}

Relation \eqref{O2} means that the convergence claimed in \eqref{O1} is not strong for $\{ \vr_n \}_{n > 0}$.

\end{Remark}

\begin{proof}

The proof is based on the method of convex integration. First, consider a division of the domain $\Omega$,
\[
\Omega = \cup_{i \in I} \Ov{\Omega}_i, \ \Omega_i \cap \Omega_j = \emptyset \ \mbox{for}\ i \ne j,
\]
where {$I$ is a finite index set, and} $\Omega_i$ are domains. Furthermore, we consider a sequence of endpoints $T_i$. Next, for each $\Omega_i$, fix $\vr_i > 0$ - a constant density distribution.
Similarly to \cite{FeKlKrMa}, \cite{LuXiXi}, we consider the following problem:
\begin{equation} \label{O3}
\begin{split}
\Div \vm_i &= 0,\\
\partial_t \vm_i + \Div \left( \frac{\vm_i \otimes \vm_i}{\vr_i} - \frac{1}{N} \frac{|\vm_i|^2}{\vr_i} \mathbb{I} \right) &= 0\\
\frac{1}{2} \frac{|\vm_i|^2}{\vr_i} &= \Lambda - p(\vr_i)\frac{N}{2}
\end{split}
\end{equation}
in $\Omega_i$, where $\Lambda > 0$ is a certain positive constant to be determined below.
The apparently overdetermined problem \eqref{O3} is supplemented by the initial--end state condition
\begin{equation} \label{O4}
\vm_i(0, \cdot) = \vm_i(T, \cdot) = 0.
\end{equation}
In addition, we impose the ``no flux'' boundary conditions specified in the weak sense as follows:
We suppose that
\begin{align} \label{O5a}
{\int_0^{T}} \int_{\Omega_i} \vm_i \cdot \Grad \varphi \ \dx \dt &= 0 
\end{align}
for any $\varphi \in C^1([0,T] \times \Ov{\Omega}_i)$ and
\begin{align}
\label{O5} {\int_0^{T}} \int_{\Omega_i} \left[ \vm_i \cdot \partial_t \bfphi + \left( \frac{\vm_i \otimes \vm_i}{\vr_i}
- \frac{1}{N} \frac{|\vm_i|^2}{\vr_i} \mathbb{I} \right): \Grad \bfphi \right] \dx \dt &= 0
\end{align}
for any $\bfphi \in C^1([0,T] \times \Ov{\Omega}_i; R^N)$.
Accordingly, solutions defined on $\Omega_i$ can be ``pasted'' together to produce a weak solution defined on the whole set $\Omega$.
Indeed, for $[\vr_i,\vm_i]$ satisfying \eqref{O4}, we can set
\begin{equation} \label{O6}
\vr = \sum_{i} 1_{\Omega_i} \vr_i,\ { \vm (t, \cdot) = \sum_{i} 1_{\Omega_i} \vm_i(t - mT)  ,\ t \in [mT,(m+1) T)
\ \mbox{for}\ m=0,1,\dots }
\end{equation}

It is a routine matter to check that $[\vr,\vm]$ defined through \eqref{O6} is a weak solution of the Euler system \eqref{i6}
{for $t \in (0, \infty)$}, satisfying the impermeability condition \eqref{i4}. Note that the momentum equation reads
\[
{\int_0^\infty}  \intO{ \left[ \vm \cdot \partial_t \bfphi + \frac{ \vm \otimes \vm}{\vr} : \Grad \bfphi + p(\vr) \Div \bfphi -
\Lambda \Div \bfphi \right] } \dt = 0,
\]
where
\[
{\int_0^\infty} \intO{ \Lambda \Div \bfphi } = 0 \ \mbox{whenever}\ \bfphi \in C^1_c({[0,\infty)} \times \Ov{\Omega}; R^d),\
\bfphi \cdot \vc{n}|_{\Omega} = 0.
\]
{Note that, in contrast with \eqref{O5} where no boundary conditions are imposed on test functions, we have
effectively used the fact $\bfphi \cdot \vc{n}|_{\partial \Omega} = 0$ here.}

Now we claim that problem \eqref{O3}--\eqref{O5} admits, in fact, infinitely many solutions as soon as
\[
0 < \underline{\vr} \leq \vr_i \leq \Ov{\vr},\ i = 1,2, \dots
\]
for certain $\Lambda = \Lambda(\underline{\vr}, \Ov{\vr}) > 0$. {Indeed we refer e.g. to Chiodaroli \cite{Chiod} or \cite{LuXiXi} for the proof.}

Finally, we consider an oscillating sequence
\[
\vr_n = \vr^n_i \in \Omega^n_i, \ \vr_n \to \vr_0 \ \mbox{weakly-(*) in}\ L^\infty(\Omega)
\ \mbox{but not strongly in}\ L^1(\Omega),
\]
with {the family of times $T_n = \frac{1}{2^n}$}, and
\[
\vm_n \ \mbox{defined on}\ {[0,\infty), \ \vm_n \left( \frac{m}{2^n}, \cdot \right) = 0, \ m=0,1,\dots.}
\]
It can be checked that $[\vr_n,\vm_n]$ enjoys the properties claimed in the conclusion of Proposition
\ref{PO1}. {Indeed we have}
\[
{\vm_n \to \vm \ \mbox{in}\ C_{\rm weak}([0,T]; L^2(\Omega; R^N)) \ \mbox{and weakly-(*) in}\ L^\infty((0,T) \times \Omega; R^N)}
\]
{for any $T > 0$. Moreover, thanks to the pointwise convergence in $L^2(\Omega; R^N)-\mbox{weak}$ at \emph{any}\ $t \geq 0$, we have}
\[
{ \vm\left( \frac{m}{2^n}, \cdot \right) = 0 \ \mbox{for any}\ m =0,1,\dots, \ n = 1,2,\dots;
\ \mbox{whence}\ \vm \equiv 0. }
\]
\end{proof}

Apparently, the limit quantity $\vr = \vr_0(x)$, $\vm \equiv 0$ is a (weak) solution of the Euler system only if $\vr_0 = \Ov{\vr}$
- a (positive) constant. Otherwise, the weak closure takes us out of the set of weak solutions. This indicates that a possibly larger class of solutions is necessary to characterize the weak closure. These are the measure--valued solutions discussed in the next section.

\section{Measure--valued solutions}
\label{M}

The concept of \emph{measure--valued} solution was introduced to capture the two major stumbling blocks to strong stability of the Euler system: (i) oscillations discussed in the previous section, and (ii) concentrations due to the kinetic energy ``blow up''.
These two phenomena are conveniently captured by the oscillation--concentration defect measure introduced by Alibert and Bouchitt\'e
\cite{AliBou}.

Gwiazda et al. \cite{GSWW} used the approach of \cite{AliBou} for the compressible Euler system. This technique requires a certain structure
of the nonlinearities to define their recession functions. This structure enforces the introduction new state variables:
the density $\vr$ and the ``weighted velocity'' $\sqrt{\vr} \vu$. It is interesting to note that similar choice of variables has been
use by Chen and Glimm \cite{CheGli} in a different context. Within this framework, Gwiazda et al. established the existence as well as the
weak--strong uniqueness principle.

The approach of \cite{GSWW} was highly simplified in \cite{FGSWW1} in the context of the compressible Navier--Stokes. The Alibert--Bouchitt\'e
defect measures have been replaced by a combination of the standard Young measure acting on the natural variables $\vr$ and $\vu$ and the concentration defect measures balanced by their dissipation counterpart in the energy inequality. This technique has been adapted
by Basari\v c \cite{Basa} to the compressible Euler system \eqref{i1} posed on a general, possibly unbounded, domain.

{Another simplification, using  rather the conservative variables $\vr,\vm\equiv\vr\vu$, has been introduced in \cite{BreFeiHof19} in order to construct a solution semiflow to the isentropic Euler system. We refer also  to \cite{BreFeiHof19B} for the application to the complete Euler system and to Section \ref{P} and Section \ref{S} below for further discussion of this subject. These developments led to the work \cite{FeiHof19}, where the underlying ideas for the notion of dissipative solution presented in the sequel can be found. This particularly straightforward formulation allowed to establish the following striking dichotomy: a weakly converging sequence of (weak) solutions to the isentropic Navier--Stokes system on $R^{N}$, $N=2,3,$ in the vanishing viscosity limit either {\bf (i)} converges strongly in the energy norm, or {\bf (ii)} the limit is not a weak solution of the associated Euler system, see \cite{FeiHof19}.}

\section{Dissipative solutions}
\label{D}

Motivated by the above mentioned results, we propose the concept of \emph{dissipative} solution adapted to the natural
\emph{conservative variables}: the density $\vr$ and the momentum $\vm$ in the Euler system \eqref{i6}.
They satisfy the following system of equations in the sense of distributions:
\begin{equation} \label{D1}
\begin{split}
\partial_t \vr + \Div \vm &= 0,\\
\partial_t \vm + \Div \left( \frac{\vm \otimes \vm}{\vr} \right) + \Grad p(\vr) &= - \Div \left( \mathfrak{R}_v +
\mathfrak{R}_p \mathbb{I} \right)\\
\partial_t \intO{ \left[ \frac{1}{2} \frac{|\vm|^2}{\vr} + P(\vr)  + \frac{1}{2} {\rm trace}[ \mathfrak{R}_v ]
+ \frac{1}{\gamma - 1} \mathfrak{R}_p \right] } &\leq 0,
\end{split}
\end{equation}
where $\mathfrak{R}_v \in L^\infty(0,T; \mathcal{M}^+(\Ov{\Omega}; R^{d \times d}_{\rm sym}))$,
$\mathfrak{R}_p \in L^\infty(0,T; \mathcal{M}^+ (\Ov{\Omega}))$ are the turbulent defect measure associated to the
convective term and the pressure, respectively. Here, the symbol $\mathcal{M}^+(\Ov{\Omega})$ denotes the space of non--negative Borel measures on $\Ov{\Omega}$, while $\mathcal{M}^+(\Ov{\Omega}; R^{d \times d}_{\rm sym})$ is the space of matrix valued (signed) measures
on $\Ov{\Omega}$ ranging in positive semi--definite matrices, meaning
\[
\mathfrak{R}_v : (\xi \otimes \xi) \in \mathcal{M}^+(\Ov{\Omega}) \ \mbox{for any}\ \xi \in R^d.
\]

Observe that dissipative solutions are weakly continuous in time, specifically,
\[
\vr \in C_{{\rm weak}}([0,T]; L^\gamma (\Omega)),\
\vm \in C_{\rm weak}([0,T]; L^{\frac{2 \gamma}{\gamma + 1}}(\Omega; R^d))
\]
so that one can correctly define the initial conditions. The boundary condition \eqref{i4} is satisfied in the weak sense
through suitable choice of the test functions in the weak formulation. The exact definition reads as follows:

\begin{Definition} \label{DD1}

We say that
\[
\vr \in C_{\rm weak}([0,T]; L^\gamma (\Omega)),\ \vr \geq 0,\ \vm \in C_{\rm weak}([0,T]; L^{\frac{2 \gamma}{\gamma + 1}}
(\Omega; R^d)),
\]
is a \emph{dissipative solution} to the Euler system \eqref{i1}--\eqref{i5} if there exist
turbulent defect measures
\[
\mathfrak{R}_v \in L^\infty(0,T; \mathcal{M}^+(\Ov{\Omega}; R^{d \times d}_{\rm sym})),\
\mathfrak{R}_p \in L^\infty(0,T; \mathcal{M}^+(\Ov{\Omega}))
\]
such that the following holds:
\begin{equation} \label{D2}
\left[ \intO{ \vr \varphi } \right]_{t = 0}^{t = \tau} =
\int_0^\tau \intO{ \Big[ \vr \partial_t \varphi + \vm \cdot \Grad \varphi \Big] } \dt
\end{equation}
for any $0 < \tau < T$, and any $\varphi \in C^1_c([0,T) \times \Ov{\Omega})$;
\begin{equation} \label{D3}
\begin{split}
\left[ \intO{ \vm \cdot \bfphi  } \right]_{t = 0}^{t = \tau} &=
\int_0^\tau \intO{ \left[ \vm \cdot \partial_t \bfphi + \left( \frac{\vm \otimes \vm}{\vr} : \Grad \bfphi \right) + p(\vr)
\Div \bfphi \right] } \dt \\ &+
\int_0^\tau \int_{\Ov{\Omega}}\Grad \bfphi : \D \Big[ \mathfrak{R}_v + \mathfrak{R}_p \mathbb{I} \Big]  \dt
\end{split}
\end{equation}
for any $0 < \tau < T$, and any $\varphi \in C^1_c([0,T) \times \Ov{\Omega}; R^d)$, $\bfphi \cdot \vc{n}|_{\partial \Omega} = 0$;
\begin{equation} \label{D4}
\begin{split}
&\left[ \psi \left( \intO{ \left[ \frac{1}{2} \frac{|\vm|^2}{\vr} + P(\vr) \right] } +
\int_{\Ov{\Omega}} \frac{1}{2} \D \, {\rm trace}[\mathfrak{R}_v] + \int_{\Ov{\Omega}} \frac{1}{\gamma - 1} \D \mathfrak{R}_p
\right) \right]_{t = \tau_1 -}^{t = \tau_2 +} \\ &\leq \int_{\tau_1}^{\tau_2}
\partial_t \psi \left( \intO{ \left[ \frac{1}{2} \frac{|\vm|^2}{\vr} + P(\vr) \right] } +
\int_{\Ov{\Omega}} \frac{1}{2} \D \, {\rm trace}[\mathfrak{R}_v] + \int_{\Ov{\Omega}} \frac{1}{\gamma - 1} \D \mathfrak{R}_p
\right) \dt
\end{split}
\end{equation}
for any $0 \leq \tau_1 \leq \tau_2 < T$, and any $\psi \in C^1_c[0,T)$, $\psi \geq 0$.
\end{Definition}

\begin{Remark} \label{DR1}

In \eqref{D4}, the initial value of the energy is set
\[
\intO{ \left[ \frac{ |\vm_0|^2 }{\vr_0} + P(\vr_0) \right] }.
\]

\end{Remark}

Although the system \eqref{D1} is apparently underdetermined due to the presence of the turbulent defect measures, it reduces to
\eqref{i6}, \eqref{i7}, meaning $\mathfrak{R}_v = \mathfrak{R}_p = 0$ as soon as $\vr$ and $\vm$ are continuously differentiable
and $\vr \geq \underline{\vr} > 0$ is bounded below away from zero. Indeed we can introduce the velocity $\vu = \frac{1}{\vr} \vm
\in C^1$, whereas the continuity equation is satisfied in the classical sense:
\[
\partial_t \vr + \Div (\vr \vu) = 0.
\]
Next, as $\vu$ can be used as a test function in the momentum equation \eqref{D3}, we easily deduce
\begin{align} \label{D4a}
\begin{aligned}
\intO{ \left[ \frac{1}{2} \frac{|\vm|^2}{\vr} + P(\vr) \right] (\tau, \cdot) } &=
\intO{ \left[ \frac{1}{2} \frac{|\vm_0|^2}{\vr_0} + P(\vr_0) \right] }\\& +
\int_{0}^\tau \int_{\Ov{\Omega}} \Grad \vu : \D \left[ \mathfrak{R}_v + \mathfrak{R}_p \mathbb{I} \right] \dt.
\end{aligned}
\end{align}
This expression may be subtracted from the energy inequality \eqref{D4} to obtain
\begin{equation} \label{D5}
\int_{\Ov{\Omega}} \left[ \frac{1}{2} \D \, {\rm trace}[\mathfrak{R}_v] + \int_{\Ov{\Omega}} \frac{1}{\gamma - 1} \D \mathfrak{R}_p
\right](\tau) \leq - \int_{0}^\tau \int_{\Ov{\Omega}} \Grad \vu : \D \left[ \mathfrak{R}_v + \mathfrak{R}_p \mathbb{I} \right] \dt
\end{equation}
Thus a direct application of Gronwall's lemma yields the desired conclusion $\mathfrak{R}_v = \mathfrak{R}_p = 0$.
We have shown the following result.

\begin{Theorem} \label{DT1}

Let $\Omega \subset R^d$ be a bounded domain of class $C^1$. Suppose that a dissipative solution $\vr$, $\vm$ is continuously differentiable in $[0,T) \times \Ov{\Omega}$ and $\vr \geq \underline{\vr} > 0$.

Then $\mathfrak{R}_v = \mathfrak{R}_p = 0$ and $\vr$, $\vm$ is a classical solution of the Euler system.

\end{Theorem}

A short inspection of \eqref{D5} shows that $C^1$ regularity is not really necessary. In fact, it is enough that
the symmetric velocity gradient
\[
\mathbb{D} \vu \equiv \frac{ \Grad \vu + \Grad \vu^t }{2}
\]
satisfies a one sided Lipschitz condition, specifically,
\begin{equation} \label{D6}
\mathbb{D} \vu + d \mathbb{I} \geq 0 \ \mbox{for certain}\ d \in L^1(0,T).
\end{equation}
Indeed, as $\mathfrak{R}_v + \mathfrak{R}_p \mathbb{I}$ is positively definite, we get
\[
\begin{split}
- \int_{0}^\tau \int_{\Ov{\Omega}} &\Grad \vu : \D \left[ \mathfrak{R}_v + \mathfrak{R}_p \mathbb{I} \right] \dt
= - \frac{1}{2} \int_{0}^\tau \int_{\Ov{\Omega}} \mathbb{D} \vu : \D \left[ \mathfrak{R}_v + \mathfrak{R}_p \mathbb{I} \right] \dt
\\ &= - \frac{1}{2} \int_{0}^\tau \int_{\Ov{\Omega}} (\mathbb{D} \vu + d \mathbb{I}) : \D \left[ \mathfrak{R}_v + \mathfrak{R}_p \mathbb{I} \right] \dt + \frac{1}{2} \int_0^\tau d\int_{\Ov{\Omega}}  \D \, {\rm trace}\left[ \mathfrak{R}_v + \mathfrak{R}_p \mathbb{I} \right] \dt \\
&\leq \frac{1}{2} \int_0^\tau d\int_{\Ov{\Omega}}  \D \, {\rm trace}\left[ \mathfrak{R}_v + \mathfrak{R}_p \mathbb{I} \right] \dt.
\end{split}
\]
Consequently, validity of Theorem \ref{DT1} can be extended to the class of dissipative solutions satisfying the energy
equality \eqref{D4a} together with the one--sided Lipschitz condition \eqref{D6} for the velocity field. Sufficient conditions for validity of the energy equality of the compressible Euler system have been studied in \cite{FeGwGwWi} in the case of periodic boundary conditions. It turns out that
\eqref{D4a} remains valid if $\vr$, $\vm$, and $\vu$ enjoy extra Besov--type regularity, specifically:
\begin{equation} \label{D7}
\begin{split}
\vr &\in L^\infty((0,T) \times \Omega),\ \vr \geq \underline{\vr} > 0,\
\vm \in L^\infty((0,T) \times \Omega; R^d),\\
\vr, \ \vm, \ \vu \equiv \frac{\vm}{\vr} &\in B^{\alpha, \infty}_3 ((0,T) \times \Omega; R^d), \ \alpha > \frac{1}{3},\\
\vr, \ \vm &\in L^\infty(0,T; B^{\beta, \infty}_q (\Omega; R^d)) \ \mbox{for some}\ \beta > 0, \ q > 1.
\end{split}
\end{equation}

\begin{Remark} \label{DR2}

The symbol $B^{\alpha, \infty}_q(Q)$ denotes the Besov space endowed with the norm
\[
\| v \|_{B^{\alpha, \infty}_q(Q)}
= \| v \|_{L^q(Q)} + \sup_{\xi \in Q} \frac{ \| v(\cdot + \xi) - v(\cdot) \|_{L^q(Q \cap (Q - \xi)) }}{|\xi|^\alpha}.
\]

\end{Remark}

Now, observe that the impermeability condition \eqref{i4}, if imposed on the cube
\[
\Omega = (-1,1)^d,
\]
can be transformed to the periodic boundary conditions working with classes of functions with certain symmetry, see Ebin \cite{EB},
and \cite{FMPS}. Summing up the previous observations, we obtain the following extension of Theorem \ref{DT1}.

\begin{Theorem} \label{DT2}

Let
\[
\Omega = (-1,1)^d
\]
be the cube. Suppose that $[\vr,\vm]$ is a dissipative solution of the Euler system belonging to the class \eqref{D7}. In addition,
let the velocity $\vu$ satisfy the one--sided Lipschitz condition
\[
\int_0^T \intO{\Big( - \xi \cdot \vu \, (\xi \cdot \Grad) \varphi + d |\xi|^2 \varphi \Big)} \dt \geq 0,\
d \in L^1(0,T),
\]
for any $\xi \in R^d$, $\varphi \in C^1_c((0,T) \times \Omega)$.

Then $\mathfrak{R}_v = \mathfrak{R}_p = 0$ and, consequently, $[\vr,\vm]$ is a weak solution of the Euler system.

\end{Theorem}

\section{Construction of dissipative solutions}
\label{C}

Dissipative solutions appear as a limit of various approximation schemes. To simplify presentation, we consider the
\emph{periodic} boundary condition, meaning the spatial domain $\Omega$ is identified with the flat torus
\[
\Omega = \left( [-1,1]|_{\{ - 1,1 \}} \right)^d.
\]

The approximate solutions typically solve a system of equations:
\begin{equation} \label{C1}
- \intO{ \vr_{0,n} \varphi } =
\int_0^T \intO{ \Big[ \vr_n \partial_t \varphi + \vm_n \cdot \Grad \varphi \Big] } \dt + E_{1,n}[\varphi]
\end{equation}
for any $\varphi \in C^1_c([0,T) \times {\Omega})$;
\begin{align}
- \intO{ \vm_{0,n} \cdot \bfphi  }  &=
\int_0^T \intO{ \left[ \vm_n \cdot \partial_t \bfphi + \left( \frac{\vm_n \otimes \vm_n}{\vr_n} : \Grad \bfphi \right) + p(\vr_n)
\Div \bfphi \right] } \dt \nonumber\\ &+ E_{2,n}[\bfphi]\label{C2}
\end{align}
for any $\bfphi \in C^1_c([0,T) \times {\Omega}; R^N)$;
\begin{equation} \label{C3}
\begin{split}
&- \intO{ \left[ \frac{1}{2} \frac{|\vm_{0,n}|^2}{\vr_{0,n}} + P(\vr_{0,n}) \right] } \leq \int_{0}^{T}
\partial_t \psi \intO{ \left[ \frac{1}{2} \frac{|\vm_n|^2}{\vr_n} + P(\vr_n) \right] } +
E_{3,n}[\psi]
\end{split}
\end{equation}
for any $\psi \in C^1_c[0,T)$, $\psi \geq 0$, $\psi(0) = 1$.

The terms $E_{1,n}$, $E_{2,n}$, $E_{3,n}$ represent \emph{consistency errors}. Furthermore, we suppose
\begin{equation} \label{C4}
E_{1,n}[\varphi] \to 0,\ E_{2,n}[\bfphi] \to 0,\
E_{3,n}[\psi] \to 0 \ \mbox{as}\ n \to \infty
\ \mbox{for fixed}\ [\varphi, \bfphi, \psi].
\end{equation}
Moreover, we require that
\begin{equation} \label{C5}
E_{3,n}[\psi] \aleq c\left(\| \psi \|_{L^\infty(0,T)} \right) \ \mbox{uniformly for} \ n \to \infty.
\end{equation}

The approximate solutions $[\vr_n, \vm_n]$ can be obtained via a numerical scheme or a suitable physically relevant approximation.
We may consider a \emph{viscosity approximation}:
\begin{equation} \label{C6}
\begin{split}
\partial_t \vr_n + \Div \vm_n = 0, \ \vr_n(0, \cdot) &= \vr_{n,0}, \\
\partial_t \vm_n + \Div \left( \frac{ \vm_n \otimes \vm_n }{\vr_n} \right) + \Grad p(\vr_n) &= \frac{1}{n} \Div \mathbb{S}_n,\
\vm_n(0, \cdot) = \vm_{0,n},
\end{split}
\end{equation}
together with the relevant energy balance
\begin{equation} \label{C7}
\frac{{\rm d}}{{\rm d}t} \intO{ \left[ \frac{1}{2} \frac{|\vm_n|^2 }{\vr_n} + P(\vr_n) \right] }
+ \frac{1}{n} \intO{ \mathbb{S}_n: \mathbb{D} \vu_n } \leq 0,
\end{equation}
where the velocity field $\vu_n$ satisfies $\vr_n \vu_n = \vm_n$. We suppose the viscous stress depends in a \emph{monotone}
way on the velocity gradient $\mathbb{D}$, meaning
\[
\mathbb{S}_n : \mathbb{D} \vu_n = F(\mathbb{D} \vu_n) + F^* (\mathbb{S}_n),
\]
where $F$ is a convex l.s.c. function on $R^{d \times d}_{\rm sym}$ and $F^*$ its conjugate. If, for instance, ${\rm Dom}[F] =
R^{d \times d}_{\rm sym}$, $F(0) = 0$, $F \geq 0$, the conjugate $F^*$ is non--negative and superlinear. Accordingly, we may set
\[
E_{1,n} = E_{3,n} = 0,\ E_{2,n}[\bfphi] = \frac{1}{n} \left| \int_0^T \intO{ \mathbb{S}_n : \Grad \bfphi } \dt \right|,
\]
whereas the desired estimates follow from the energy balance \eqref{C7}.

Our next goal is to perform the limit for $n \to \infty$ in \eqref{C1}--\eqref{C3}. The key tool is the energy inequality \eqref{C3}
yielding, together with the consistency bound \eqref{C5}, the uniform bounds
\begin{equation} \label{C8}
\intO{ \left[ \frac{1}{2} \frac{|\vm_n|^2}{\vr_n} + P(\vr_n) \right] }
\leq c({\rm data}) \ \mbox{uniformly for}\ t \in (0,T) \ \mbox{and}\ n \to \infty,
\end{equation}
in particular,
\begin{align} \label{C9}
\begin{aligned}
\vr_n &\to \vr \ \mbox{weakly-(*) in}\ L^\infty(0,T; L^\gamma (\Omega)),\\
\vm_n &\to \vm \ \mbox{weakly-(*) in}\ L^\infty(0,T; L^{\frac{2 \gamma}{\gamma + 1}}(\Omega; R^d)),
\end{aligned}
\end{align}
for suitable subsequences as the case may be. Note that the function
\[
[\vr, \vm] \in [0, \infty) \times R^d \mapsto \frac{|\vm|^2}{\vr} =
\left\{ \begin{array}{l} 0 \ \mbox{whenever}\ \vm = 0, \\
\frac{|\vm|^2}{\vr} \ \mbox{for}\ \vr > 0, \\
\infty \ \mbox{otherwise} \end{array} \right.
\]
is a convex l.s.c. function.

Next, we have, again for a subsequence,
\[
p(\vr_n) \to \Ov{p(\vr)} \ \mbox{weakly-(*) in}\ L^\infty(0,T; \mathcal{M}^+(\Omega)).
\]
Moreover, as $p$ is convex, we have
\[
0 \leq p(\vr) \leq \Ov{p(\vr)}, \ \mathfrak{R}_p \equiv \Ov{p(\vr)} - p(\vr)
\in L^\infty(0,T; \mathcal{M}^+(\Omega)).
\]

By the same token,
\[
\frac{\vm_n \otimes \vm_n}{\vr_n} \to \Ov{ \frac{\vm \otimes \vm}{\vr}}
\ \mbox{weakly-(*) in}\ L^\infty(0,T; \mathcal{M}(\Omega; R^{d \times d}_{\rm sym})),
\]
and we set
\[
\mathcal{R}_v \equiv \Ov{ \frac{\vm \otimes \vm}{\vr}} -
\frac{\vm \otimes \vm}{\vr}.
\]
The crucial observation now is that
\[
\begin{split}
\mathcal{R}_v : (\xi \otimes \xi) &= \lim_{n \to \infty} \left[ \frac{ \vm_n \otimes \vm_n }{\vr_n} : (\xi \otimes \xi) \right] -
\frac{\vm \otimes \vm}{\vr} : (\xi \otimes \xi)\\ &=
\lim_{n \to \infty} \left[ \frac{|\vm_n \cdot \xi|^2}{\vr_n} - \frac{|\vm \cdot \xi|^2}{\vr} \right] \geq 0
\end{split}
\]
due to convexity. We therefore conclude that
\[
\mathfrak{R}_v \in L^\infty(0,T; \mathcal{M}^+(\Omega; R^{d \times d}_{\rm sym})).
\]

Finally, it is a routine matter to show that the limit $[\vr,\vm]$ satisfies \eqref{D2}--\eqref{D4}, meaning, it is a dissipative solution of the Euler system.

\section{Properties of the solution set}
\label{P}

Dissipative solutions are not uniquely determined by the initial data unless they enjoy certain extra regularity property similar to \eqref{D7}. However,
we report the following \emph{weak--strong} uniqueness principle proved in \cite[Theorem 2.1]{FeiGhoJan}.

\begin{Theorem} \label{PT1}
Let
\[
\Omega = \left( [-1,1]|_{\{ - 1,1 \}} \right)^d
\]
be the flat torus. Suppose that the Euler system \eqref{i1}--\eqref{i3} admits a weak solution $\widetilde{\vr}$, $\widetilde{\vm}=
\widetilde{\vr} \vU$ belonging to the regularity class:
\begin{equation} \label{P1}
\begin{split}
\widetilde{\vr} &\in B^{\alpha, \infty}_{p}((\delta,T) \times \Omega)) \cap C([0,T]; L^1(\Omega)),\\
\vU &\in B^{\alpha, \infty}_{p}((\delta,T) \times \Omega; R^d)) \cap C([0,T]; L^1(\Omega; R^d)),
\ \mbox{for any}\ \delta > 0,
\end{split}
\end{equation}
with
\[
\alpha > \frac{1}{2}, \ p \geq \frac{4 \gamma}{\gamma - 1};
\]
\[
0<
\underline{\vr} \leq \widetilde{\vr}(t,x) \leq \Ov{\vr}, \ | \vU(t,x)| \leq \Ov{U}\ \mbox{for a.a.} \ (t,x) \in (0,T) \times \Omega;
\]
\[
\intO{ \left[ - \xi \cdot \vU(\tau, \cdot) (\xi \cdot \Grad) \varphi  + d(\tau) |\xi|^2 \varphi \right] } \geq 0,\ d \in L^1(0,T),
\]
for any $\xi \in R^d$, and any $\varphi \in C(\Omega)$, $\varphi \geq 0$. Let $\vr$, $\vm$ be a dissipative solution starting from the initial data
\[
\vr(0, \cdot) = \widetilde{\vr}(0,\cdot),\ \vm(0, \cdot) = \widetilde{\vm}(0, \cdot).
\]

Then $\mathfrak{R}_p = \mathfrak{R}_v = 0$, and $\vr = \widetilde{\vr}$, $\vm = \widetilde{\vm}$.

\end{Theorem}

In the remaining part of this section, we examine the properties of the solution set for fixed finite energy initial data:
\begin{equation} \label{P2}
\vr_0 \in L^\gamma(\Omega),\ \vm_0 \in L^{\frac{2 \gamma}{\gamma + 1}}(\Omega; R^d),\ \intO{ \left[ \frac{1}{2} \frac{ |\vm_0|^2 }{\vr_0} + P(\vr_0) \right] } \leq E_0 < \infty.
\end{equation}
Let
\[
\begin{split}
\mathcal{U} [\vr_0, \vm_0] = &\Big\{ [\vr, \vm] \ \Big| \ \vr \in C_{\rm weak}([0,T]; L^\gamma(\Omega)),\
\vm \in C_{\rm weak}([0,T]; L^{\frac{2 \gamma}{\gamma + 1}}(\Omega; R^d)) \\ &\mbox{is a dissipative solution of the Euler system},
\ \vr(0, \cdot) = \vr_0, \ \vm(0, \cdot) = \vm_0 \Big\}
\end{split}
\]
be the set of all dissipative solutions in the sense of Definition \ref{DD1} starting from the initial data $[\vr_0, \vm_0]$.

We claim that for any $[\vr_0, \vm_0]$ satisfying \eqref{P2}:

\begin{itemize}

\item
$\mathcal{U}[\vr_0, \vm_0]$ is non--empty;
\item
$\mathcal{U}[\vr_0, \vm_0]$ is convex;
\item
$\mathcal{U}[\vr_0, \vm_0]$ is compact with respect to the metric topology on bounded sets in
\[
\Big[
C_{\rm weak}([0,T]; L^\gamma(\Omega)) \cap L^\infty(0,T; L^\gamma(\Omega)) \Big] \times
\Big[
C_{\rm weak}([0,T]; L^{\frac{2 \gamma}{\gamma + 1}}(\Omega)) \cap L^\infty(0,T; L^{\frac{2 \gamma}{\gamma + 1}}(\Omega)) \Big].
\]

\end{itemize}

The fact that the solution set is non--empty was proved in Section \ref{C}. Compactness can be shown by the same arguments as
the proof of existence. Finally, as a convex combination of two (non--negative) measures is a measure, it is easy to check that the set
$\mathcal{U}[\vr_0, \vm_0]$ is convex.

\section{Selection criteria, admissible solutions}
\label{S}

As we have observed in the previous section, the set of dissipative solutions $\mathcal{U}[\vr_0, \vm_0]$ emanating from the
initial data $[\vr_0, \vm_0]$ is non--empty, convex, and compact with respect to the weak topology on the trajectory space.
Unfortunately, there are numerous examples furnished by the method of convex integration showing the set $\mathcal{U}[\vr_0, \vm_0]$
is not a singleton.

Several criteria could be proposed to rule out the irrelevant solutions. We discuss shortly the \emph{maximal dissipation principle} asserting that the physical solution dissipates the (mechanical) energy at maximal rate. Given $[\vr_1, \vm_1]$, $[\vr_2, \vm_2]$
we define a relation
\[
[\vr_1, \vm_1] \prec [\vr_2, \vm_2]
\]
if
\[
\frac{1}{2} \frac{|\vm_1|^2}{\vr_1} + P(\vr_1) + \frac{1}{2} {\rm trace}[\mathfrak{R}^1_v] +
\frac{1}{\gamma - 1} \mathfrak{R}^1_p \leq
\frac{1}{2} \frac{|\vm_2|^2}{\vr_2} + P(\vr_2) + \frac{1}{2} {\rm trace}[\mathfrak{R}^2_v] +
\frac{1}{\gamma - 1} \mathfrak{R}^2_p
\]
in the sense of measures on $[0,T] \times \Ov{\Omega}$.

\begin{Definition} \label{SD2}

Let the initial data $[\vr_0, \vm_0]$ be given. We say that a dissipative solution $[\vr, \vm]$ is \emph{admissible} if it
is minimal with respect to the relation $\prec$. Specifically, if $[\widetilde{\vr}, \widetilde{\vm}]$ is another
dissipative solution starting from the same initial data such that
\[
[\widetilde{\vr}, \widetilde{\vm}] \prec [\vr, \vm],
\]
then
\[
\frac{1}{2} \frac{|\vm|^2}{\vr} + P(\vr) + \frac{1}{2} {\rm trace}[\mathfrak{R}_v] +
\frac{1}{\gamma - 1} \mathfrak{R}_p =
\frac{1}{2} \frac{|\widetilde{\vm}|^2}{\widetilde{\vr}} + P(\widetilde{\vr}) + \frac{1}{2} {\rm trace}[\widetilde{\mathfrak{R}}_v] +
\frac{1}{\gamma - 1} \widetilde{\mathfrak{R}}_p
\]
in $[0,T] \times \Ov{\Omega}$.

\end{Definition}

It is easy to see that an admissible solution always exist. It is enough to minimize the energy functional
\[
\int_0^T \intO{ \left[ \frac{1}{2} \frac{|\vm|^2}{\vr} + P(\vr) \right] } \dt +
\int_0^T \int_{\Ov{\Omega}} \left[ \frac{1}{2}\D \, \mbox{trace}[\mathfrak{R}_v] + \frac{1}{\gamma - 1} \D\,
\mathfrak{R}_p \right]
\]
over the set of all dissipative solutions $[\vr, \vm]$ with the associated turbulence defects $\mathcal{R}_v$, $\mathcal{R}_p$
starting from the initial data $[\vr_0, \vm_0]$.

{Finally, we point out that a suitable choice of a family of \emph{selection criteria} gives rise to a
\emph{semiflow selection} and conditional well posedness. The basic idea goes back to the Krylov \cite{KrylNV}, where a general selection procedure has been proposed in the context of Markov semigroups. Similar approach in the deterministic setting was used by
Cardona and Kapitanskii \cite{CorKap}. Subsequently, the method was adapted to the compressible Euler system in \cite{BreFeiHof19}.
More precisely, the state variables being enhanced by the associated energy $E$, there is a measurable mapping
$$
U:[t,\vr_{0},\vm_{0},E_{0}]\mapsto [\vr(t),\vm(t),E(t)],
$$
such that $[\vr,\vm,E]$ solves the Euler system (in the sense of dissipative solutions) with the initial data given by $[\vr_{0},\vm_{0},E_{0}]$ and the semigroup property
$$
U[t_{1}+t_{2},\vr_{0},\vm_{0},E_{0}]=U[t_{2},U[t_{1},\vr_{0},\vm_{0},E_{0}]]\ \text{for any}\  t_{1},t_{2}\geq 0,
$$
holds true.
The interested reader may consult \cite{BreFeiHof19} for details.}

\def\cprime{$'$} \def\ocirc#1{\ifmmode\setbox0=\hbox{$#1$}\dimen0=\ht0
  \advance\dimen0 by1pt\rlap{\hbox to\wd0{\hss\raise\dimen0
  \hbox{\hskip.2em$\scriptscriptstyle\circ$}\hss}}#1\else {\accent"17 #1}\fi}


\begin{thebibliography}{10}

\bibitem{AliBou}
J.~J. Alibert and G.~Bouchitt\'e.
\newblock Non-uniform integrability and generalized {Y}oung measures.
\newblock {\em J. Convex Anal.}, {\bf 4}(1):129--147, 1997.

\bibitem{Basa}
D.~Basari{\' c}.
\newblock Vanishing viscosity limit for the compressible {N}avier--{S}tokes
  system via measure-valued solutions.
\newblock {\em Arxive Preprint Series}, {\bf arXiv 1903.05886}, 2019.

\bibitem{BreFeiHof19B}
D.~Breit, E.~Feireisl, and M.~Hofmanov{\' a}.
\newblock Dissipative solutions and semiflow selection for the complete {E}uler
  system.
\newblock {\em Arxive Preprint Series}, {\bf arXiv 1904.00622}, 2019.

\bibitem{BreFeiHof19}
D.~Breit, E.~Feireisl, and M.~Hofmanov{\' a}.
\newblock Solution semiflow to the isentropic {E}uler system.
\newblock {\em Arxive Preprint Series}, {\bf arXiv 1901.04798}, 2019.

\bibitem{CorKap}
J.E. Cardona and L.~Kapitanskii.
\newblock Semiflow selection and {M}arkov selection theorems.
\newblock {\em Arxive Preprint Series}, {\bf arXiv 1707.04778v1}, 2017.

\bibitem{CheGli}
G.~Q. Chen and J.~Glimm.
\newblock {K}olmogorov-type theory of compressible turbulence and inviscid
  limit of the {N}avier--{S}tokes equations in ${R}^3$.
\newblock {\em Arxive Preprint Series}, {\bf arXiv 1809.09490}, 2018.

\bibitem{Chiod}
E.~Chiodaroli.
\newblock A counterexample to well-posedness of entropy solutions to the
  compressible {E}uler system.
\newblock {\em J. Hyperbolic Differ. Equ.}, {\bf 11}(3):493--519, 2014.

\bibitem{ChiDelKre}
E.~Chiodaroli, C.~{D}e {L}ellis, and O.~Kreml.
\newblock Global ill-posedness of the isentropic system of gas dynamics.
\newblock {\em Comm. Pure Appl. Math.}, {\bf 68}(7):1157--1190, 2015.

\bibitem{ChiKre}
E.~Chiodaroli and O.~Kreml.
\newblock On the energy dissipation rate of solutions to the compressible
  isentropic {E}uler system.
\newblock {\em Arch. Ration. Mech. Anal.}, {\bf 214}(3):1019--1049, 2014.

\bibitem{ChKrMaSwI}
E.~Chiodaroli, O.~Kreml, V.~M{\' a}cha, and S.~Schwarzacher.
\newblock Non–uniqueness of admissible weak solutions to the compressible
  {E}uler equations with smooth initial data.
\newblock {\em Arxive Preprint Series}, {\bf arXiv 1812.09917v1}, 2019.

\bibitem{DelSze3}
C.~De~Lellis and L.~Sz{\'e}kelyhidi, Jr.
\newblock On admissibility criteria for weak solutions of the {E}uler
  equations.
\newblock {\em Arch. Ration. Mech. Anal.}, {\bf 195}(1):225--260, 2010.

\bibitem{EB}
D.~B. Ebin.
\newblock Viscous fluids in a domain with frictionless boundary.
\newblock {\em Global Analysis - Analysis on Manifolds, H. Kurke, J. Mecke, H.
  Triebel, R. Thiele Editors, Teubner-Texte zur Mathematik 57, Teubner,
  Leipzig}, pages 93--110, 1983.

\bibitem{FeiGhoJan}
E.~Feireisl, S.~S. Ghoshal, and A.~Jana.
\newblock On uniqueness of dissipative solutions to the isentropic {E}uler
  system.
\newblock {\em Arxive Preprint Series}, {\bf arXiv 1903.11687}, 2019.

\bibitem{FGSWW1}
E.~Feireisl, P.~Gwiazda, A.~{\'S}wierczewska-Gwiazda, and E.~Wiedemann.
\newblock Dissipative measure-valued solutions to the compressible
  {N}avier--{S}tokes system.
\newblock {\em Calc. Var. Partial Differential Equations}, 55(6):55:141, 2016.

\bibitem{FeGwGwWi}
E.~Feireisl, P.~Gwiazda, A.~\'Swierczewska-Gwiazda, and E.~Wiedemann.
\newblock Regularity and energy conservation for the compressible {E}uler
  equations.
\newblock {\em Arch. Ration. Mech. Anal.}, {\bf 223}(3):1375--1395, 2017.

\bibitem{FeiHof19}
E.~Feireisl and M.~Hofmanov{\' a}.
\newblock On the vanishing viscosity limit of the isentropic {N}avier--{S}tokes
  system.
\newblock {\em Arxive Preprint Series}, {\bf arXiv 1905.02548}, 2019.

\bibitem{FeKlKrMa}
E.~Feireisl, C.~Klingenberg, O.~Kreml, and S.~Markfelder.
\newblock On oscillatory solutions to the complete {E}uler system.
\newblock {\em Arxive Preprint Series}, {\bf arXiv 1710.10918}, 2017.

\bibitem{FMPS}
E.~Feireisl, {\v S}.~Matu{\v s}{\accent23u}-{N}e{\v c}asov{\' a},
  H.~Petzeltov{\' a}, and I.~Stra{\v s}kraba.
\newblock On the motion of a viscous compressible flow driven by a
  time-periodic external flow.
\newblock {\em Arch. Rational Mech. Anal.}, {\bf 149}:69--96, 1999.

\bibitem{GSWW}
P.~Gwiazda, A.~\'Swierczewska-Gwiazda, and E.~Wiedemann.
\newblock Weak-strong uniqueness for measure-valued solutions of some
  compressible fluid models.
\newblock {\em Nonlinearity}, {\bf 28}(11):3873--3890, 2015.

\bibitem{KrylNV}
N.~V. Krylov.
\newblock The selection of a {M}arkov process from a {M}arkov system of
  processes, and the construction of quasidiffusion processes.
\newblock {\em Izv. Akad. Nauk SSSR Ser. Mat.}, {\bf 37}:691--708, 1973.

\bibitem{LuXiXi}
T.~Luo, C.~Xie, and Z.~Xin.
\newblock Non-uniqueness of admissible weak solutions to compressible {E}uler
  systems with source terms.
\newblock {\em Adv. Math.}, {\bf 291}:542--583, 2016.

\bibitem{SMO}
J.~Smoller.
\newblock {\em Shock waves and reaction-diffusion equations}.
\newblock Springer-Verlag, New York, 1967.

\bibitem{TAN}
A.~Tani.
\newblock On the first initial-boundary value problem of compressible viscous
  fluid motion.
\newblock {\em Publ. RIMS Kyoto Univ.}, {\bf 13}:193--253, 1977.

\end{thebibliography}

\end{document}